\documentclass[11pt]{article}
\usepackage{paper,asb}

\usepackage{amsmath}
\usepackage{amsfonts}
\usepackage{amssymb}
\usepackage{amsthm}
\usepackage{hyperref}
\usepackage{subcaption}
\usepackage{cite}

\usepackage{upgreek} 
\usepackage{bm}

\numberwithin{equation}{section}

\usepackage{booktabs}
\usepackage{caption}
\usepackage{float}
\usepackage{titlesec}
\usepackage{capt-of}
\usepackage{arydshln}

\author{Jiahao Shi \and James C. Spall}
\date{}

\newcommand{\papertitle}{Difference Between Cyclic and Distributed Approach in Stochastic Optimization for Multi-agent System}
\newcommand{\paperauthora}{Jiahao Shi}
\newcommand{\paperauthoraaffiliation}{University of Michigan}
\newcommand{\paperauthoraemail}{jiahaos@umich.edu}
\newcommand{\paperauthorb}{James C. Spall}
\newcommand{\paperauthorbaffiliation}{Johns Hopkins University}
\newcommand{\paperauthorbemail}{james.spall@jhuapl.edu }

\title{\papertitle}

\author{\paperauthora\footnotemark[2] \footnotemark[1]
   \and \paperauthorb\footnotemark[3]}

\begin{document}

\maketitle

\renewcommand{\thefootnote}{\fnsymbol{footnote}}
\footnotetext[2]{\paperauthoraaffiliation. (\url{\paperauthoraemail})}
\footnotetext[3]{\paperauthorbaffiliation. (\url{\paperauthorbemail})}
\footnotetext[1]{Corresponding author.}
\renewcommand{\thefootnote}{\arabic{footnote}}

\begin{abstract}{
Many stochastic optimization problems in multi-agent systems can be decomposed into smaller subproblems or reduced decision subspaces. The cyclic and distributed approaches are two widely used strategies for solving such problems. In this manuscript, we review four existing methods for addressing these problems and compare them based on their 
suitable problem frameworks and update rules. 
}
\end{abstract}

\maketitle




\section{Introduction}
\label{sec:intro}

Stochastic optimization refers to a class of problems where 
only noisy observations of the objective function $L(\bm\uptheta)$ or its gradient may be obtained. Formally, the problem is 
\begin{align}
    \min_{\bm{\uptheta} \in \mathbb{R}^p} L(\bm{\uptheta}), 
\label{eq:1}
\end{align}
 where $L:\mathbb{R}^p \to \mathbb{R}$ is a differentiable loss function with its gradient denoted as $\bm{g}(\cdot)$.  
 Many algorithms are developed to solve such problems \cite{robbins1951stochastic,duchi2011adaptive,kingma2014adam,johnson2013accelerating,ghadimi2013stochastic,kiefer1952stochastic,nelder1965simplex,powell2002uobyqa,larson2016stochastic,chen2018stochastic,spall1992multivariate,berahas2019derivative,cao2023first,larson2024novel,sun2023trust,shi2021sqp}. Among these algorithms, \cite{larson2016stochastic,chen2018stochastic,spall1992multivariate,berahas2019derivative,cao2023first,larson2024novel,sun2023trust} are specifically designed for derivative-free setting where only the noisy observations of the objective function $L(\bm\uptheta)$ are available. 
In some applications, such as multi-agent systems \cite{olfati2007consensus,dorri2018multi}, each agent typically has incomplete information about the environment, making general stochastic optimization algorithms inappropriate. To deal with such problems,  the cyclic approach and the  distributed approach have been considered over
the last years.

The cyclic (or seesaw, alternating, block coordinate)  approach to optimization
 divides the full
parameter vector into multiple subvectors and
the process proceeds by sequentially optimizing the criterion
of interest with respect to each of the subvectors while
holding the other subvectors fixed. The decision variable
$\bm{\uptheta} \in \mathbb{R}^p$ is  represented here as 
\begin{align}
    \bm{\uptheta} =  [\bm{t}_{1}^T, \bm{t}_{2}^T, \cdots, \bm{t}_{A}^T]^T, 
\label{eq:division}
\end{align}
where $\bm{t}_{i}^T \in \mathbb{R}^{p_i}$ is a subvector with $p_i$ being its  number of components in $\bm{\uptheta}$ satisfying $\sum_{i=1}^A p_i = p$, and $A$ is the total number of subvectors. At each iteration, only one subvector is updated. 
Several studies have examined cyclic optimization (not necessarily applied to the multi-agent systems) in deterministic settings \cite{bertsekasandtsitsiklis1989,tseng2001,nesterov2012efficiency, spallcyclicseesaw2012,wright2015}. 
As surveyed in \cite{wright2015}, the pattern for updating subvectors of the cyclic approach  can be either deterministic or random. The version of random subvector selection can be considered as a special case of stochastic approximation algorithm due to its randomness in updating the subvector at each iteration. However, this algorithm is not suitable to the stochastic optimization problem (\ref{eq:1}) in which case we can only obtain noisy function measurements of $L(\bm \uptheta)$. Although  the research work \cite{liu2014asynchronous} seemingly considers the cyclic approach in the stochastic setting, the randomness  is actually a consequence of  the update pattern. As a result, it does not apply the cyclic approach to the stochastic optimization problem of noisy loss or gradient measurements. References \cite{tsitsiklis1994asynchronous,borkar1998asynchronous} provide significant contributions to the cyclic algorithm in stochastic approximation. However, both works impose relatively strong assumptions on the noise (e.g., bounded variance) and do not offer theoretical results such as asymptotic normality or relative efficiency compared to their non-cyclic counterparts. Moreover, these approaches do not support blockwise coordinate updates or randomized update patterns. To address these limitations, the generalized cyclic stochastic approximation (GCSA) algorithm proposed in  \cite{cuevas2017cyclic,hernandez2014cyclic} explores the cyclic approach in stochastic optimization, assuming the existence of a unique solution $\bm{\uptheta}^*$ for (\ref{eq:1}).


 Cyclic procedures have been applied in various fields, including control problems in robotics \cite{canutescu2003}, computer vision \cite{leenpark2008}, compressed sensing \cite{liandosher2009}, parameter estimation \cite{haaland2010statistical}, and target location estimation \cite{liandpetriouli2014}. In this manuscript, we focus on multi-agent problems where each agent, whose decision is represented by a subvector of $\bm \uptheta$, faces uncertainty about the entire system, as reflected in noisy measurements of $L$. We refer to this problem framework as the \textit{cyclic framework}. 
 To find the minimizer of $L$ in this setting, each agent controls and updates its local decision variable. 
 An example of such a problem, along with the application of the algorithm, will be discussed later.

However, the term multi-agent problem can be ambiguous, as it is also used to describe problems formalized as 
\begin{align}
    \min_{\bm{\uptheta} \in \mathbb{R}^p} L(\bm{\uptheta}) = \sum_{i=1}^A L_i(\bm \uptheta).
    \label{distribute} 
\end{align}
This popular framework  is based on a sum of agent contributions,  where in a network of $A$ agents, each agent is endowed with a local loss function $L_i:\mathbb{R}^p \to \mathbb{R}$ known only to agent $i$. It is often assumed that the each $L_i$ is local strictly convex \cite{nedic2009distributed,berahas2018balancing}. We refer to this problem framework as the \textit{distributed framework}. 
There is a growing literature on developing 
algorithms for solving problem (\ref{distribute}) both in the deterministic setting  \cite{bertsekasandtsitsiklis1989,nedic2009distributed,boyd2011distributed,peng2016arock,shi2015extra,nedic2017achieving,berahas2018balancing,berahas2023balancing} and the stochastic setting  \cite{tsitsiklis1986distributed,ram2010distributed,ram2009incremental,ouyang2013stochastic,yang2016parallel,pu2021distributed}. In the stochastic setting, it is commonly assumed that each agent $i$ has access only to a noisy estimate of its local function or its gradient. 
Existing methods solving such problems 
include the  Distributed Stochastic Approximation/Subgradient (DSA) algorithm \cite{ram2010distributed}, the Incremental Stochastic Approximation/Subgradient (ISA) algorithm \cite{ram2009incremental}, the  Stochastic Alternating Direction Method of Multipliers (SADMM) algorithm \cite{ouyang2013stochastic}, and  the Stochastic Parallel Successive Convex Approximation-based (SPSCA) algorithm \cite{yang2016parallel}. The SADMM algorithm is a general extension of the well-known Alternating Direction Method of Multipliers (ADMM) algorithm \cite{boyd2011distributed} in the stochastic setting and considers a constrained stochastic optimization problem. Although the cyclic approach is  used in the SADMM algorithm in the sense where two variables of the augmented  Lagrangian function are updated separately, we consider the SADMM algorithm as a dual method instead of a cyclic method. 
The SPSCA algorithm uses a decomposition approach and updates the local estimate of the global decision variable by solving a sequence of strongly convex subproblems.
The subproblem that each agent aims to solve is based on a surrogate function of the local loss function  $L_i(\cdot)$.  
This paper provides a detailed review of both the DSA and ISA algorithms. 
The distributed approach has been applied to specific areas including 
control systems \cite{gazi2011swarm,cao2012overview,shamma2007cooperative}, wireless sensor networks \cite{bullo2018lectures,derakhshan2019review}, smart grid \cite{pipattanasomporn2009multi},  
power engineering \cite{mcarthur2007multi,mcarthur2007multi2},  autonomous vehicle systems \cite{gazi2011swarm,peterson2014simulation,botts2016multi,botts2020three,chandler2001uav,dong2018theory}, and   machine learning \cite{tsianos2012consensus,wooldridge2009introduction,forero2010consensus}. By using multi-agent systems, certain global objectives can be achieved through sensing, exchange of information using
communication, computation and control of local objectives.


Throughout the paper, we consider multi-agent optimization problems, where each agent might be a sensor, a robot, a vehicle, or a computing node, as those problems involve multiple agents cooperatively minimizing a loss function. In this sense, we cover a wide range of real-world problems in the literature  \cite{gazi2011swarm,cao2012overview,shamma2007cooperative,tsianos2012consensus,wooldridge2009introduction,derakhshan2019review,pipattanasomporn2009multi,mcarthur2007multi,mcarthur2007multi2,beard2002coordinated,chandler2001uav} including both \textit{cyclic framework} and \textit{distributed framework}. In many applications, agents can cooperate in centralized or decentralized settings. A centralized approach is based on the assumption that a central station is available and powerful enough to control a whole group of agents to produce a globally optimal solution. By contrast, a decentralized approach allows each agent to be optimized separately in which each agent communicates and computes without a central station so that it produces a locally optimal solution for each agent \cite{peterson2014simulation,cao2012overview}. Compared with the centralized control algorithms, decentralized control algorithms are more robust to accidental agent failures and communication link disruptions \cite{tsitsiklis1984problems}. This paper focuses on the decentralized setting 
because of the benefits mentioned above.

The main contribution of this paper is to compare the difference between the cyclic and distributed approach in stochastic multi-agent problems. 
The cyclic approach can be applied to general stochastic optimization problem \eqref{eq:1} when the parameter vector can be divided into subvectors. The distributed approach, however, is mostly applied in problem \eqref{distribute} under a network of multiple agents since the loss function has to be a sum of agent contributions.  
In multi-agent system,  both the cyclic and distributed approaches assume that each agent has incomplete knowledge about the environment in the decentralized setting, but are different in which part of information available to each agent. Consequently, the update rule of these two algorithms are different. 

The outline of the paper is as follows. In Section~\ref{sec:algorithm}, we introduce the algorithms to be compared. Section~\ref{sec:diff} provides a comparison of these algorithms. Finally, concluding remarks are presented in Section~\ref{sec:conclusion}.

\section{Overview of Existing Methods}
\label{sec:algorithm}
Stochastic approximation (SA) is a powerful class of iterative algorithms for minimizing a loss function, $L(\bm \uptheta)$, with respect to a parameter vector $\bm \uptheta$, when only noisy observations of $L(\bm \uptheta)$ or its gradient are available. This section provides an overview of four existing stochastic approximation algorithms. These methods share some algorithmic similarities and can all be applied in multi-agent systems. To assist practitioners in distinguishing between related algorithms, we review them in this section and highlight their differences in Section~\ref{sec:diff}.



\subsection{Generalized Cyclic Stochastic Approximation (GCSA) Algorithm}
\label{sec:GCSA}
In this section, we present the GCSA algorithm that is studied in \cite{cuevas2017cyclic}, which considers the cyclic approach and its generalization in classical SA theory with the theoretical results derived in \cite{cuevas2017cyclic}. The proposed GCSA algorithm is fully described in Algorithm \ref{alg.GCSA},  
 which sequentially optimizes a subvector of $\bm \uptheta \in \mathbb{R}^p$ while holding the other elements fixed. Let the set $\mathcal{S}\equiv \{1,\dots,p\}$ be divided into $A$ (not necessarily disjoint) subsets $\{\mathcal{S}_i\}_{i=1}^A$  satisfying  $\bigcup_{i=1}^A\mathcal{S}_i=\mathcal{S}$. This definition coincides with the  block coordinate approach defined in  \cite{tseng2001}.  By this definition, $\mathcal{S}_i$ contains the coordinates in the $i$th subvector. Note that $\{\mathcal{S}_i\}_{i=1}^A$ are not necessarily disjoint in the GCSA algorithm.  
For simplicity, we only present a special case of the GCSA algorithm in this paper and assume that $\bigcap_{i=1}^A\mathcal{S}_i=\emptyset$. We do this simplification since this assumption is often attainable in multi-agent optimization systems in which each agent has its own decision variables. This special case updates the subvectors in a strictly cyclic way (from agent $1, 2, \cdots$ to $A$).
Note that the GCSA algorithm (Algorithm 1 in \cite{cuevas2017cyclic}) is more general in that it allows both the subvector and the times of each subvector to be updated  to be randomized per iteration.  
   
For $i=1, \ldots, A$, let $\bm{g}^{(i)} \in \mathbb{R}^p$ denote the information of $\bm{g}$ for the $i$th subvector, the generally non-zero components of $\bm{g}^{(i)}$ correspond to the $i$th subvector. In particular, let $g_{j}^{(i)}, g_{j}$ denote the $j$th entry of the vector $\bm{g}^{(i)}, \bm{g}$ respectively, and let $g_{j}^{(i)}=g_{j}$ if $j \in S_{i}$ and $g_{j}^{(i)}=0$ otherwise.  For example, suppose that $\bm{g} = [g_1, g_2,\cdots,g_p]^T$  and $S_1 = \{1,2\}$, then $\bm{g}^{(1)} = [g_1, g_2,0,\cdots,0]^T$.  With this notation, it is easy to observe that the recursive form (line 4 in Algorithm~\ref{alg.GCSA}) only updates a subvector of $\hat{\bm{\uptheta}}_k$ with $\hat{\bm g}^{(i)}(\cdot)$ being the gradient estimate of subvector $i$ of
the true gradient ${\bm g}(\cdot)$ with zeros elsewhere. Note that we allow a modest assumption for the gradient estimate that fits into those assumptions needed in classical stochastic approximation (SA) theory. See, for example,  \cite[Sect. 4.3.2]{spall2005introduction} or \cite[Sect. 2.2]{kushner2012stochastic}. As a result, the gradient estimate can be stochastic gradient (SG), finite difference SA (FDSA), or simultaneous perturbation SA (SPSA).   Under certain assumptions, it is proved in \cite{cuevas2017cyclic} and \cite{hernandez2014cyclic} that for Algorithm~\ref{alg.GCSA}, we have $\hat{\bm{\uptheta}}_k \to {\bm{\uptheta}}^*$ w.p.1.  


\begin{algorithm}[ht]
  \caption{Strictly Cyclic Version of GCSA}
  \label{alg.GCSA}  
  \begin{algorithmic}[1]
    \Require $\hat{\bm{\uptheta}}_0 \in \mathbb{R}^p$, $\{a_k\}_{k\geq 0}$. Set $k=0$. 
    \While{stopping criterion has not been reached}
    \State{Set $\hat{\bm\uptheta}_{k0}=\hat{\bm{\uptheta}}_k$;}
    \For{$i=1,\dots,A$}
    \State{Update $\hat{\bm\uptheta}_{ki}= \hat{\bm\uptheta}_{k,i-1}- {a}_{k}\ \hat{\bm{g}}^{(i)}(\hat{\bm\uptheta}_{k,i-1})$.} 
    \EndFor
    \State{Let $\hat{\bm{\uptheta}}_{k+1}\equiv\hat{\bm{\uptheta}}_{kA}$;
	 }
	\State{Set $k=k+1$.}
    \EndWhile
      \end{algorithmic}
\end{algorithm}

\subsection{Distributed Stochastic Approximation (DSA) Algorithm}
\label{sec:DSA}
In this section, we review the distributed approach of stochastic approximation. In the literature, a series of distributed stochastic approximation (DSA) algorithms apply the SA method to the stochastic multi-agent optimization problem (\ref{distribute}) including \cite{ram2010distributed,sahu2018distributed,sergeenko2020advanced,bianchi2013performance}. By using different methods to obtain a gradient estimate $\hat{\bm g}(\cdot)$ of $L(\cdot)$,  \cite{ram2010distributed} studies a distributed stochastic subgradient algorithm,  \cite{sahu2018distributed} studies a distributed finite difference SA algorithm, and  \cite{sergeenko2020advanced} studies a distributed simultaneous perturbation SA algorithm.

The stochastic multi-agent optimization problem (\ref{distribute}) assumes that the function $L_i(\cdot)$ is known only partially to agent $i$ in the sense that the agent can only obtain a noisy gradient estimate of its own contribution to the loss function defined as $\hat{\bm g_i}(\cdot)$.
Agents are allowed to exchange information with their peers via a network. 
The communication network is often determined by directed or undirected graphs. For simplicity, this paper only considers the network determined by undirected graphs. 
In an undirected graph, an edge $(i,j)$ denotes that agent $j$ and agent $i$ can communicate with each other. Through communication, agent $j$ can obtain information of agent $i$ and vice versa. 

Determined by whether the topology of the graph is time-invariant or dynamic, the communication network can be static or dynamic \cite{nedic2009distributed,cao2012overview,shamma2007cooperative}. Determined by a dynamic network, the scalars $\bm W_{ij}(k) \in [0,1]$ are defined as the non-negative weights satisfying  $\sum_{j=1}^A \bm W_{ij}(k) =1$ that agent $i$ assigns to  agent $j$ in the $k$th iteration. If the communication network is static, then we have $\bm W_{ij}(k) = \bm W_{ij}$ for all $k$. Let $N_i(k)$ denote the set of neighboring agents of agent $i$ whose current iterates are available to agent $i$ in the $k$th iteration, then the DSA algorithm can then be described as Algorithm \ref{DGD}. 

To cooperatively optimize the global decision variable  $\bm \uptheta \in \mathbb{R}^p$, each agent has its local estimate  $\hat{\bm{\uptheta}}_{ki}\in \mathbb{R}^p$ initialized at  $\hat{\bm{\uptheta}}_{0i}$ for $i=1,\cdots,A$ ($A$ is the total number of agents). The general idea of how the DSA algorithm works is that each agent updates  $\hat{\bm{\uptheta}}_{ki}$ using only the information of its local loss function. The agent also communicates with other agents to share the local information, which helps ensure   that all the local estimates are converging to the global optimal $\bm \uptheta^*$.  In each iteration, each agent first communicates with its neighboring agents and obtain the weighted averaged information to achieve consensus (line 2 of Algorithm \ref{DGD}). Then, it updates its local estimate using the weighted averaged information and the gradient estimate $\hat{\bm g_i}(\cdot)$ (line 3 of Algorithm \ref{DGD}). By performing this type of local communication and computation steps, the consensus in the estimates is achieved under certain assumptions (see \cite{ram2010distributed}), satisfying 
\begin{align}
    \lim_{k \to \infty} \| \hat{\bm{\uptheta}}_{ki} - \hat{\bm{\uptheta}}_{k} \| = 0 \text{ for all } i = 1,\cdots,A, 
\end{align}
where the convergence is with probability 1 and in mean square. In addition,  we have the global convergence result 
\begin{align}
    \lim_{k \to \infty} \| \hat{\bm{\uptheta}}_{k} - \bm{\uptheta}^* \| = 0
\end{align}
with probability 1. 


\begin{algorithm}[ht]                     
  \caption{Distributed Stochastic Approximation (DSA) Algorithm}      
  \label{DGD}                
  \begin{algorithmic}[1]  
    \Require $\hat{\bm{\uptheta}}_{0i}\in \mathbb{R}^p$ for $i=1,\cdots,A$, $\{a_k\}_{k\geq 0}$. Set $k=0$.
    \While{stopping criterion has not been reached, each agent $i$}
    \State Define $\hat{\bm{\upeta}}_{ki} = \sum_{j \in N_i(k)} \bm W_{ij}(k) \hat{\bm{\uptheta}}_{ki}$;
    \State Update $\hat{\bm{\uptheta}}_{k+1,i} = \hat{\bm{\upeta}}_{ki} - a_k \hat{\bm g_i}( \hat{\bm{\upeta}}_{ki})$;
    \State Set $k=k+1$.
    \EndWhile
    \State Define $\hat{\bm{\uptheta}}_{k} = \frac{1}{A} \sum_{i=1}^A  \hat{\bm{\uptheta}}_{ki}$;
  \end{algorithmic}
\end{algorithm}


\subsection{Other Related Algorithms}
In this paper, we want to compare the cyclic and the distributed approach in solving stochastic optimization problems. The GCSA algorithm and the DSA algorithm are different in both the problem  framework and the update rule. We will discuss about those differences in detail in Section~\ref{sec:diff}. Before making a direct comparison, we present two algorithms that fall between the GCSA algorithm and the DSA algorithm. The purpose of this section is to include other methods in the literature that  resemble the two approaches of above, but which differ in important, but subtle, ways. 

\subsubsection{DSA--Specialization (DSA--S) Algorithm}
\label{sec:DSAS}

A special case of the DSA algorithm is considered in \cite[Example \uppercase\expandafter{\romannumeral3}]{tsitsiklis1986distributed}. 
In the context of multi-agent system, this special case does not assume that each agent has a local loss function $L_i(\bm{\uptheta})$ as in problem (\ref{distribute}), but assumes that each agent specializes a single component of the global decision variable $\bm{\uptheta}$. 
For simplicity, let the $i$th component of $\bm{\uptheta}$ be specialized to agent $i$. Let $\hat{\bm g}^{(i)}(\bm{\uptheta})$ be the $i$th component of $\hat{\bm g}(\bm{\uptheta})$. Notice that the notation of $\hat{\bm g}^{(i)}(\bm{\uptheta})$ here coincides with $\hat{\bm g}^{(i)}(\bm{\uptheta})$ in the GCSA algorithm in the case that each subvector has only one component of $\bm{\uptheta}$.  Following \cite{tsitsiklis1986distributed}, we call this special case the DSA--Specialization (DSA--S) algorithm, and present it as Algorithm~\ref{DSAspecial}. Similar to the DSA algorithm, each agent has its local estimate  $\hat{\bm{\uptheta}}_{ki}\in \mathbb{R}^p$ initialized at  $\hat{\bm{\uptheta}}_{0i}$ for $i=1,\cdots,A$. 
In the computation step, agent $i$ only update its particular component in $\hat{\bm{\uptheta}}_{ki}$. The other components of $\hat{\bm{\uptheta}}_{ki}$ are only updated in the communication step. 

Both Algorithms~\ref{DGD} and \ref{DSAspecial} can be implemented in parallel. This is because at the $k$th iteration, the update of each agent only relies on the information associated with $\hat{\bm{\uptheta}}_{ki}$ for $i=1,\cdots,A$, suggesting that the all agents can communicate or compute simultaneously. The main difference between Algorithms \ref{DGD} and \ref{DSAspecial} is that when minimizing the global loss $L(\bm{\uptheta})$, Algorithm \ref{DGD} distributes the local loss  $L_i(\bm{\uptheta})$ among   agents while Algorithm~\ref{DSAspecial} distributes  the processing of the components of vector   $\bm{\uptheta}$ among   agents. 


\begin{algorithm}[ht]
  \caption{DSA--Specialization (DSA--S) Algorithm}
  \label{DSAspecial}
  \begin{algorithmic}[1]
    \Require $\hat{\bm{\uptheta}}_{0i} \in \mathbb{R}^p$ for $i=1,\cdots,A$, $\{a_k\}_{k\geq 0}$. Set $k=0$.
    \While{stopping criterion has not been reached, each agent $i$}
      \State Define $\hat{\bm{\upeta}}_{ki} = \sum_{j \in N_i(k)} \bm W_{ij}(k) \hat{\bm{\uptheta}}_{ki}$;
      \State Update $\hat{\bm{\uptheta}}_{k+1,i} = \hat{\bm{\upeta}}_{ki} - a_k \hat{\bm g}^{(i)}( \hat{\bm{\upeta}}_{ki})$;
      \State Set $k=k+1$.
    \EndWhile
  \end{algorithmic}
\end{algorithm}

\subsubsection{Incremental Stochastic Approximation Algorithm} 
\label{sec:ISA}
Like the DSA algorithm, the incremental stochastic approximation (ISA) algorithm  is proposed to solve the finite-sum stochastic optimization problem  (\ref{distribute}). By using different methods to obtain a gradient estimate $\hat{\bm g}(\cdot)$ of $L(\cdot)$,  \cite{ram2009incremental}  studies an incremental stochastic subgradient algorithm and  \cite{yagan2006distributed} studies an incremental simultaneous perturbation SA algorithm. 
The main difference between the ISA and the DSA algorithm is that the ISA algorithm passes the estimate from one agent to another and only one agent updates at any given time. By contrast, all the agents can update simultaneously in the DSA algorithm. 

We present a special case of the ISA algorithm, cyclic incremental stochastic approximation (CISA) algorithm  in Algorithm \ref{CISS}, where the update rule follows a strictly cyclic pattern (from agent $1, 2, \cdots$ to $A$) at each iteration. There is another version of the ISA algorithm, called Markov randomized incremental stochastic approximation (MRISA), that allows the update rule to be randomized. The MRISA algorithm assumes that agents can only communicate  with their neighbors and assuming agent $i$ communicates with its neighboring $j$ with probability $Q_{ij}(k)$. Here the probability $Q_{ij}(k)$ is similar to the weight $W_{ij}(k)$ defined in \eqref{sec:DSA}.  For the ISA algorithm, under some conditions, we have the convergence result 
\begin{align}
    \lim_{k \to \infty} \hat{\bm{\uptheta}}_{k} = \bm{\uptheta}^*
\label{convergece:CISS}    
\end{align}
with probability 1.

Compared to the DSA algorithm, the ISA algorithm can be applicable to the problems where not every weight in the network is known but there exists some rule to update the agent sequentially. For instance, in the CISA algorithm, the minimum requirement for the structure of the graph is that a cyclic path exists such communication is allowed along the path.


\begin{algorithm}[ht]                     
  \caption{Cyclic Incremental Stochastic Approximation (CISA) Algorithm}    
  \label{CISS}                
  \begin{algorithmic}[1]  
    \Require $\hat{\bm{\uptheta}}_0 \in \mathbb{R}^p$, $\{a_k\}_{k\geq 0}$. Set $k=0$.
    \While{stopping criterion has not been reached}
    \State Set $\hat{\bm{\uptheta}}_{k0}=\hat{\bm{\uptheta}}_k$;
    \For{$i=1,\dots,A$}
    \State Update $\hat{\bm{\uptheta}}_{ki}= \hat{\bm{\uptheta}}_{k,i-1}- {a}_{k} \hat{\bm{g}}_{i}\left(\hat{\bm{\uptheta}}_{k,i-1}\right)$;
    \EndFor
    \State Let $\hat{\bm{\uptheta}}_{k+1} \equiv \hat{\bm{\uptheta}}_{kA}$;
    \State Set $k=k+1$;
    \EndWhile
  \end{algorithmic}
\end{algorithm}

\section{Difference Between Cyclic and Distributed Approach}
\label{sec:diff}
The similarities and differences for the algorithms  presented in Section~\ref{sec:algorithm}  may not always be evident. The goal of this section is to illustrate the differences among these algorithm  with respect to both the framework of the optimization problem and the update rule. An overview of the difference among these algorithms is shown in Table \ref{table1}. The framework of the optimization problem is different in whether to solve the problem (\ref{eq:1}) or (\ref{distribute}). The update rule is different in whether the agents update in parallel or sequential mode. Note that both the GCSA and the ISA algorithms have a generalized form of the cyclic update by allowing a randomized update pattern. For comparison purposes, only the strictly cyclic version of both algorithms is discussed. 
\begin{table}[h!] 
 \centering
 \caption{Comparison of algorithms for solving stochastic multi-agent problems.}\resizebox{\textwidth}{14mm}{
 \begin{tabular}{ccccc}
  \toprule  
  Algorithm & distributed framework & cyclic framework  & distributed update & cyclic update \\
  \midrule 
  GCSA  &   & \checkmark &  & \checkmark\\
    DSA &    \checkmark & & \checkmark & \\
      DSA-S  &   & \checkmark & \checkmark &  \\
   ISA&  \checkmark &  &  & \checkmark \\
  \bottomrule 
 \end{tabular}}
 \label{table1}
\end{table}
 
\subsection{Differences in the Framework}
\label{sec:difference_framework}
The difference in the framework of the optimization problem divides the algorithms into two classes --- \textit{cyclic framework} and \textit{distributed framework}. We categorize the
GCSA and the DSA-S algorithm as algorithms with \textit{cyclic framework}, the DSA and the ISA algorithm as algorithms with \textit{distributed framework}. The algorithms with \textit{distributed framework} are designed for solving stochastic optimization problems where the loss function (\ref{distribute}) to minimize is a sum of local loss functions only available to each agent. In  contrast,  the other class of algorithms (cyclic) do not require this structure. They are designed to minimize a more general stochastic loss function (\ref{eq:1}) without any form of local objectives.  
For example, in the context of multi-agent systems, each agent should be assumed to share the same objective instead of owning a local objective. Therefore, neither framework includes the other. It is important to notice that the algorithms with \textit{cyclic framework} are not restricted to multi-agent systems. The \textit{cyclic framework} can  be applied to general stochastic optimization problems as long as the decision variable can be divided into several subvectors. However, in the context here (per Section~\ref{sec:intro}) the \textit{distributed framework} automatically indicates a multi-agent system. 
Another distinction in the general framework is that in algorithms with a distributed structure, each agent maintains a copy of the global decision variable $\bm{\uptheta}$, serving as a local estimate $\hat{\bm{\uptheta}}_{ki} \in \mathbb{R}^p$. Unlike in \textit{distributed framework}, where each agent updates all components of the decision variable, algorithms with a \textit{cyclic framework} are designed for scenarios where each agent only has access to specific components of $\bm{\uptheta}$. The estimate for the global decision variable can be expressed as
\begin{align*}
    \hat{\bm{\uptheta}}_{k} = [\hat{\bm{t}}_{k1}^T, \hat{\bm{t}}_{k2}^T, \cdots, \hat{\bm{t}}_{kA}^T]^T,
\end{align*}
where $\hat{\bm{t}}_{ki} \in \mathbb{R}^{p_i}$ is an estimate for local decision variable of agent $i$, $p_i$ is the number of components in $\bm{\uptheta}$ possessed by agent $i$. In this case, agent $i$ can only update its local decision variable $\hat{\bm{t}}_{ki}$.  

In addition, the two classes of algorithms are different in computing local gradient estimates (or gradient approximation in the derivative-free setting). Both the gradient estimates for agent $i$, $\hat{\bm{g}}_{i}(\cdot)$,  in algorithms within the \textit{distributed framework} are estimates of the gradient of the local loss function $L_i(\cdot)$ with respect to all decision variables. However,  the gradient estimates for agent $i$, $\hat{\bm{g}}^{(i)}(\cdot)$, in algorithms within the \textit{cyclic framework} are estimates of the gradient  of the full loss function $L(\cdot)$ with respect to certain components of the global decision variables.  

Considering these differences between the algorithms with \textit{cyclic framework} and \textit{distributed framework}, different class of algorithms are used to solve different problems.
When 
each agent only has a certain contribution of the global objective and does not have the information about the objectives of its peers, the algorithms in the \textit{cyclic framework} cannot work. When the loss function is of a general form instead of a summation, or when only partial gradient information is available, the algorithms with \textit{distributed framework} cannot work. We end this section with the following two examples to show the different problems that different algorithms can solve. 


\textbf{Example 1 (Distributed Regression }\cite{ram2009incremental}\textbf{).}
Consider $A$ agents that sense a time invariant spatial field. Let $r_{ik}$ be the noisy  measurement of the true values of local objective function $R_{i}$ made by agent $i$ in time slot $k$. Let $s_i$ be the location of  agent $i$. Let $h(s; \bm{\uptheta})$ be a candidate model for the spatial field that is selected based on a priori information and parameterized by $\bm{\uptheta}$. 
The problem in regression is to estimate the parameters in the candidate model based on the collected measurements $r_{ik}$, i.e., to determine the value for $\bm{\uptheta}$ that best describes the spatial field. In a minimum mean-square sense,  the parameter value $\bm{\uptheta}^{*}$ corresponding to the best model satisfies the following relation:
\begin{align}
    \bm{\uptheta}^{*} = \underset{\bm{\uptheta} \in \mathbb{R}^p}{\operatorname{argmin}} \sum_{i=1}^{A} \mathbb{E}\left[\left(R_{i}-h\left(s_{i}, \bm{\uptheta}\right)\right)^{2}\right].
\end{align}
Let $L_i(\bm{\uptheta}) = \mathbb{E}\left[\left(R_{i}-h\left(s_{i}, \bm{\uptheta}\right)\right)^{2}\right]$, the least-square optimization problem  then fits into (\ref{distribute}), and agent $i$ can only obtain noisy measurements of $L_i(\bm{\uptheta})$ by
\begin{align}
\label{eq.least.square}
    \hat{L}_i(\bm{\uptheta}) = \frac{1}{N_i} \sum_{k=1}^{N_i} \left(r_{ik}-h\left(s_{i}, \bm{\uptheta} \right)\right)^{2},
\end{align}
where $N_i$ is the number of measurements sensed by agent $i$. This problem is also known as distributed empirical risk minimization. Because  agent $i$ cannot obtain the measurements sensed by other agents, the $\hat{L}_j(\bm{\uptheta}) (j\ne i)$ are unavailable to agent $i$. As a result,  algorithms in the \textit{cyclic framework} are not suitable for solving this problem, as they require a global objective function accessible to all agents. 

In contrast, the \textit{distributed framework} is more appropriate. The noisy local loss function in \eqref{eq.least.square} has a closed form, allowing for straightforward computation of its gradient, 
$\hat{\bm g}_i(\bm{\uptheta})$ can therefore be easily computed.  Each agent updates its own local copy $\bm{\uptheta}_i$ through a computation and communication step. By incorporating the gradient estimates into the DSA or ISA algorithms, the distributed regression problem can be effectively solved.


\textbf{Example 2 (Multi-agent Multi-target Surveillance }\cite{botts2016multi}\textbf{).} A stochastic multi-agent, multi-target surveillance problem is first considered in \cite{peterson2014simulation}, where the loss function is an information-based representation of the expected information gain resulting from agents making specific motions. The authors use cyclic stochastic optimization methods to conduct multi-agent cooperative decision making.  This framework is further extended in \cite{botts2016multi,botts2020three} to solve the same problem in 2D and 3D spaces. The application of cyclic approach in stochastic multi-agent problems can also be found in \cite{granichin2018cyclic}
to solve a parameter tracking problem. 
In \cite[Sect. 8]{cuevas2017cyclic} and \cite[Sect. 6]{zhu2020error}, the authors consider similar multi-agent problem as in \cite{peterson2014simulation,botts2016multi,botts2020three} but under a zero-communication setting. 

Consider a  representative  anti-submarine  warfare  scenario  in  which  a  fleet  of  autonomous  undersea  vehicles  (AUVs) (denoted agents) cooperatively patrol a large region to detect, track, and classify threat submarines (denoted targets).  Assume that  the total number of agents is $A$, the total number of targets is $T$, and each target updates its state according to a prescribed trajectory or motion model. Assume the states of the targets are unknown but can only be estimated by the agents. As a result, there is no optimal steady-state solution or configuration of the agents. The estimation process for each agent that updates the state estimate for each detected target is preceded by an extended Kalman filter as outlined in \cite{peterson2014simulation}.
Further, we assume that the spatial configuration of the agents and targets entirely determines which targets can be detected and which can not. The agents optimize their configuration using decentralized
motion planning. That is, each agent decides on its heading and/or vertical displacement at each time step via minimization of a loss function.

From the optimization perspective, the goal in solving a multi-agent, multi-target surveillance problem is to minimize $L_k(\bm \uptheta_k)$, where $L_k(\cdot)$ represents a time-varying stochastic loss function and $\bm \uptheta_k$ represents the global decision variables at time $k$. In particular, $\bm \uptheta_k = [t_{k1},t_{k2},\cdots, t_{kA}]^T$ with $t_{ki}$
representing 
the moving angle of  agent $i$ at time $k$. By assuming a constant speed for agents, the location of agent $i$ at any given time could be determined by its starting point and $\{t_{ki}\}_{k\ge 0}$. A  motion planning strategy for the spatial configuration of all agents can therefore be decided by $\{ \bm \uptheta_k \}_{k\ge 0} $. Note that the global loss function $L_k(\cdot)$ also depends on the state of targets, which is also a random quantity.  Formally, 
\begin{align}
 L_k(\bm{\uptheta}_k)= - \sum_{j=1}^T  \log | \bm{F}_{kj}(\bm{\uptheta}_k) |,   \label{eq:BottsLoss}
\end{align}
where $\bm{F}_{kj}(\bm{\uptheta}_k)$ is the post-action Fisher information used to quantify information about target $j$ after
new data is acquired if all agents execute
action $\bm \uptheta_k$. $L_k(\bm \uptheta_k)$ is then the negative of the total information gain aggregated over all targets as a result of all agents executing $\bm \uptheta_k$.  It is a random quantity since it depends on knowing how agents and targets move over the time step interval and which
agents (if any) will detect target $j$ after all actions are
executed at time step $k$. Because it depends on the future information of the agent and the targets, $L_k(\bm \uptheta_k)$ is not available to us, which makes most search algorithms  ineffective. Specifically, Algorithms~\ref{DGD} and \ref{CISS} 
cannot work exactly as described above in that if agents are allowed to update their local decision variables simultaneously, it will be difficult to estimate the future information of the agent and the targets without the knowledge of which agent is tracking which target. We next show that the cyclic stochastic optimization method can be used to solve this problem.

The cyclic approach assumes that each agent specializes in updating its own movement parameters, $\bm t_{ki} \in \mathbb{R}^4$,  representing angles and speeds in a 2D space, which constitute a subvector of the decision variable $\bm \uptheta_k$. 
We define 
\begin{align}
     \tilde{\bm \uptheta}_{ki} = [\tilde{\bm t}_{k1}^T,\cdots,\tilde{\bm t}_{k,i-1}^T,\bm t_{ki}^T,\tilde{\bm t}_{k,i+1}^T,\cdots,\tilde{\bm t}_{kA}^T]^T
\label{eq:subvector}
\end{align}
as the local estimate by agent $i$ of the moving angles of all agents, 
where $\tilde{\bm t}_{kl} (l \neq i)$ is the estimated moving angle of agent $l$ realized by agent $i$ communicating or predicting, it is possible to derive an expression for the loss function dependent on $\tilde{\bm \uptheta}_{ki}$. We use $ \hat{L}_{ki} (\tilde{\bm  \uptheta}_{ki})$ to  denote the derived loss function, which can be considered as a noisy measurement of the incalculable $L_k(\bm \uptheta_k)$. In \cite{botts2016multi}, $ \hat{L}_{ki} (\tilde{\bm  \uptheta}_{ki})$ is computed based on state estimation via the extended Kalman filter. 
The cyclic approach is desirable because $\tilde{\bm t}_{kl} (l \neq i)$ are fixed quantities and $\hat{L}_{ki} (\tilde{\bm \uptheta}_{ki})$ depends only on the local decision variable $\bm t_{ki}$ for agent $i$. As a result, we have
\begin{equation}
{\frac{ \partial \hat{L}_{ki}( \tilde{\bm  \uptheta}_{ki})}{\partial \tilde{\bm  \uptheta}_{ki} }} = \left[\bm 0^T, \left(\frac{ \partial \hat{L}_{ki}( \tilde{\bm  \uptheta}_{ki})}{\partial \bm t_{ki} }\right)^T ,\bm 0^T\right]^T.
\end{equation}

It is easy to observe that the gradient estimate $\partial \hat{L}_{ki}( \tilde{\bm  \uptheta}_{ki})/ \partial \tilde{\bm  \uptheta}_{ki} $ resembles the definition of $\bm v^{(l)}$ in Section~\ref{sec:difference_framework}. We can therefore apply the GCSA algorithm. That is, each agent takes turns updating its local decision variables at each iteration as follows:
\begin{align} 
     \text{For } &i=1,\cdots,A: \notag \\ 
     & \text{agent } i \text{ collects } \tilde{\bm t}_{kl} \text{ to construct } \tilde{\bm  \uptheta}_{ki} \text{ for } l \neq i; \label{eq:commun} \\
      & \hat{\bm   \uptheta}_{k,i+1} = \hat{\bm  \uptheta}_{ki} - a_k   \left. {\frac{ \partial \hat{L}_{ki}( \tilde{\bm  \uptheta}_{ki})}{\partial \tilde{\bm  \uptheta}_{ki} }} \right|_{\bm t_{ki} = {\hat{\bm t}_{ki}}} \text{ with }  \hat{\bm  \uptheta}_{k1} = \hat{\bm  \uptheta}_{k-1,A+1}. \label{eq:update}
\end{align} 
In each step, only the moving angle of a single agent is updated. Note that in this problem, the loss function and associated optimal parameter  vector $\bm  \uptheta^* $ is generally time-varying. Hence, the stepsize $a_k$ is often chosen to be a constant sequence in order to ensure that $\hat{\bm  \uptheta}_{ki} $ is able to track $\bm  \uptheta^* $.  


\subsection{Differences in the Update Rule}
\label{sec:difference_update}



In addition to the specific problems each algorithm is designed to solve, their update rules also differ.
Unlike Section~\ref{sec:difference_framework}, the difference in update rule divides the algorithms into two different classes --- \textit{cyclic update} and  \textit{distributed update}. We categorize the
GCSA and the ISA algorithm as algorithms with \textit{cyclic update}, the DSA and the DSA--S algorithm as algorithms with \textit{distributed update}. 
It is interesting to note that the two algorithms with \textit{cyclic update} are only different in the gradient estimate, and likewise for the algorithms with \textit{distributed update}. Also, both the algorithms with \textit{cyclic framework} and \textit{distributed framework} 
 contain an algorithm with \textit{cyclic update}  and an algorithm with \textit{distributed update}. That is to say, once the framework of the optimization is determined, the algorithm to use can still be different because of different update rule. 
Although it is mentioned in Section \ref{sec:difference_framework} that the algorithms with \textit{cyclic framework} are not restricted in multi-agent systems, the update rule has no differences among all problem settings. Hence we compare the update rule of the algorithms with \textit{cyclic} and \textit{distributed update} in the context of multi-agent system to simplify the comparison. 
 
\begin{figure}[ht]
\begin{subfigure}{.5\textwidth}
  \centering
  \includegraphics[width=.75\linewidth]{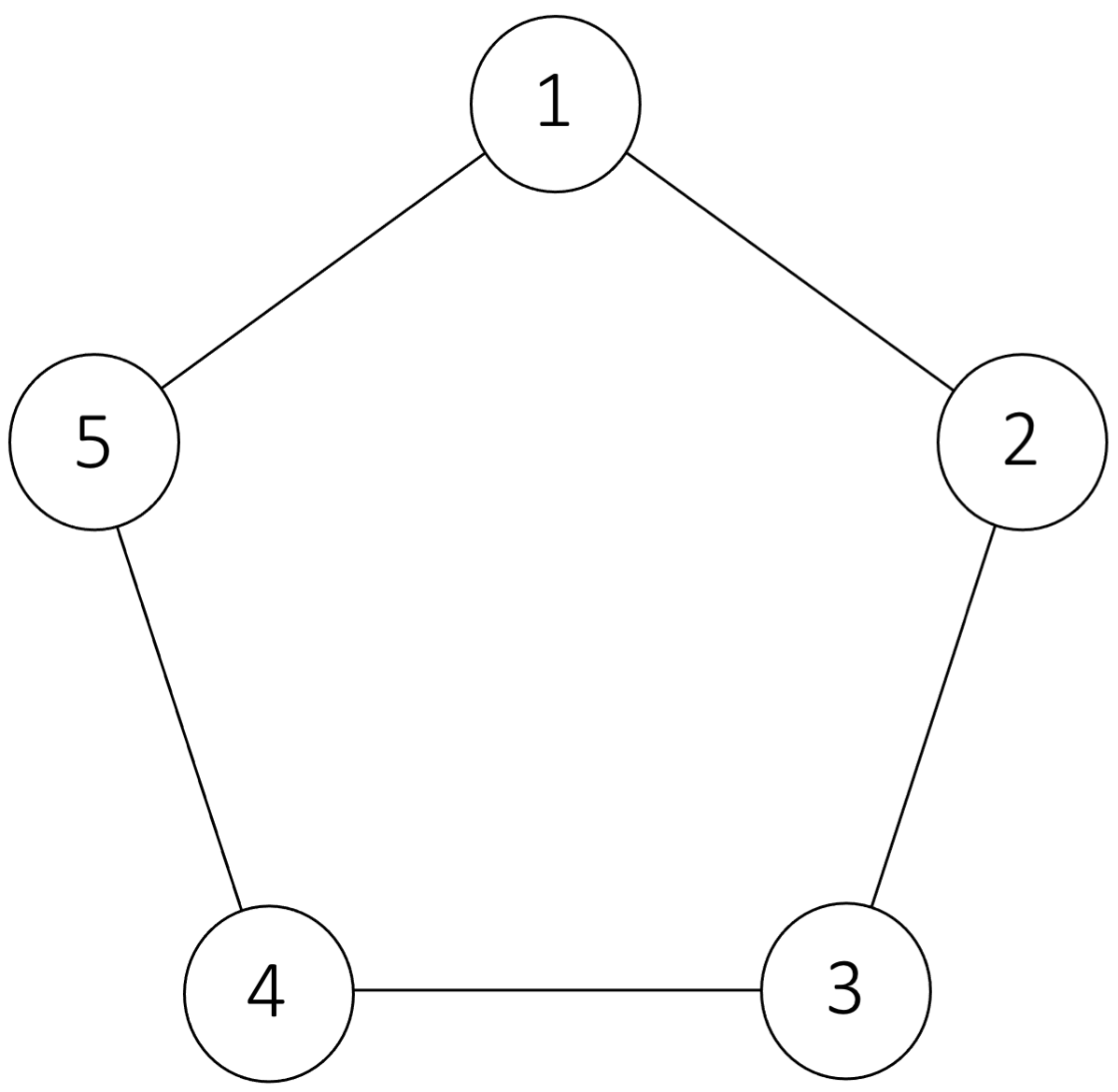}  
  \caption{}
\end{subfigure}
\begin{subfigure}{.5\textwidth}
  \centering
  \includegraphics[width=.75\linewidth]{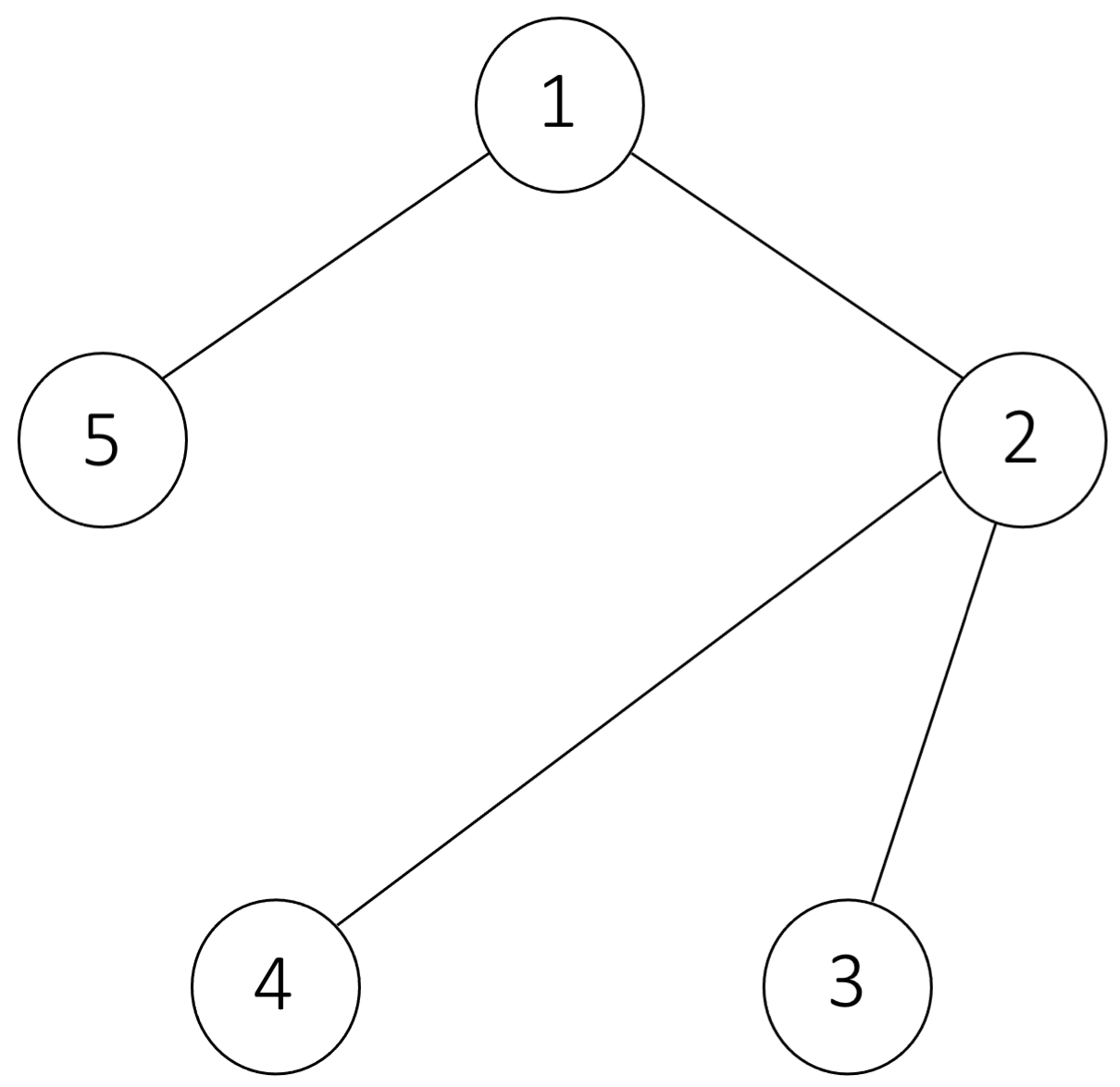}  
  \caption{}
\end{subfigure}
\caption{(a) represents a communication network with a ring structure while (b) does not  }
\label{fig:ring}
\end{figure}

Both Algorithms~\ref{alg.GCSA} and \ref{CISS} are designed for a ring (cycle) connectivity structure and this structure can make the most recent information of the peer agents available by a strictly \textit{cyclic update} pattern (one agent communicating with the next), while the algorithms with \textit{distributed update} can be applied to multi-agent problems with dynamic communication networks. Communication networks of a ring structure are  illustrated in Figure~\ref{fig:ring}. It seems that the algorithms with \textit{distributed update} are more flexible in the structure of the  communication networks. However, we argue that the algorithms with \textit{cyclic update} can also deal with  dynamic communication networks since Algorithms~\ref{alg.GCSA} and \ref{CISS} are both special cases for the GCSA and the ISA algorithm, respectively. As discussed in Section \ref{sec:GCSA} and \ref{sec:ISA}, both the GCSA and the MRISA algorithm determine which agent to update randomly at each iteration. They can therefore easily be applied to  dynamic communication networks if we know the probability of each agent communicating with its neighborhood agent at any given time.  

The distinction of algorithms with \textit{cyclic} and \textit{distributed update} is as follows. At the $k$th iteration, the intermediate estimates $\{ \hat{\bm\uptheta}_{ki} \}_{i=1}^A$ in the algorithms with \textit{cyclic update} are generated, and the update of $\hat{\bm\uptheta}_{ki}$  proceeds  only when 
$\hat{\bm\uptheta}_{k,i-1}$, the updated state of the last agent, is obtained.  Therefore, agents in these two algorithms  cannot act simultaneously and should  update sequentially. This update pattern can be categorized as the \textit{cyclic update}. By contrast, there are no intermediate estimates in the algorithms with \textit{distributed update} since the update at the $k$th iteration  only requires local estimates $\{ \hat{\bm{\uptheta}}_{ki} \}_{i=1}^A$ generated at the $(k-1)$th iteration. Because the updates only require previous information, agents in these two algorithms can act simultaneously and a parallel implementation is therefore possible. This update pattern can be categorized as the \textit{distributed update}. 


We now discuss the advantages and the disadvantages of each type of update rule (distributed or cyclic). To begin with, note that at each iteration, all of algorithms discussed in Section \ref{sec:algorithm}  require a communication step and a computation step at each iteration. The total cost is therefore produced by $A$ communication steps and $A$ computation steps. The advantage of a parallel implementation is that the communication and the computation step of each agent can be executed at the same time. As a result, implementing algorithms with \textit{distributed update} in parallel can reduce the total time cost by $1/A$ at best compared to its sequential version. 
However, this is the best conclusion we can obtain in terms of the advantage of \textit{distributed update} with  parallel implementation.  By no means does it show that algorithms with  \textit{distributed update} have a computational cost lower than algorithms with \textit{cyclic update}. Here are three reasons: 

First, the computation time  is only one component of the total time cost. Another cost is communication time. Algorithms with \textit{distributed update} require each agent to communicate with all its neighboring agents while  the algorithms with \textit{cyclic update} only require each agent to communicate with one single agent. The total communication times per iteration for the algorithms with \textit{distributed update} are therefore much greater. 

Second, even when the communication cost is cheap and the total time cost of algorithms with \textit{distributed update} can be much less per iteration, a fair comparison can be made only when we know how many iterations are required for each iteration of the algorithms  until some level of accuracy is achieved. Unfortunately, whether algorithms with \textit{distributed update} or  \textit{cyclic update} require fewer  iterations appears to be unknown in the literature. Thus it is impossible to compare the total time cost of the two classes of the algorithms.

Third, the cost of an algorithm is not fully determined by the total execution time, especially in the stochastic settings. In stochastic optimization, we often assume that the dominant cost is the noisy measurements of the loss function or its gradient \cite[p13]{spall2005introduction}.  
By this criterion, the costs of both algorithm are comparable since the number of noisy measurements of the gradient of the loss function for both algorithms per iteration is $A$.

For the reasons above, the algorithms with \textit{distributed update} are not generally the superior methods even though they can be implemented in parallel. 
We next conclude this section with a specific problem setting in which all of  algorithm discussed in Section \ref{sec:algorithm} can solve and coincide with each other. In that case, the loss function is defined in a separable form, 
\begin{align} 
  \min_{\bm{\uptheta}\in \mathbb{R}^p}  L(\bm \uptheta) = \min_{\bm{\uptheta}\in \mathbb{R}^p} \sum_{i=1}^A L_i(\bm t_i),
\label{sum3}
\end{align}
where $\bm \uptheta = [\bm t_{1}^T,\bm t_{2}^T,\cdots,\bm t_{d}^T]^T$ with each subvector $\bm t_{i}$ representing the state of agent $i$ for $i = 1, \cdots, d$. Suppose that we can only obtain the noisy measurements of $L_i(\bm t_i)$, denoted as 
$\hat{L}_i(\bm t_i)$. 
From this structure, each local loss function  only depends on its local decision variables, we then have
\begin{align}
    \bm{g}^{(i)}\left(\bm \uptheta\right) = \bm{g}_{i}\left(\bm \uptheta\right) = \left[\bm{0}^T, \left(\frac{\partial L_i(\bm t_i) }{\partial \bm t_i}\right)^T,\bm{0}^T\right]^T,
\end{align} 
where $\hat{\bm{g}}^{(i)}\left(\bm \uptheta\right)$ in algorithms with \textit{cyclic update} is the estimate of $\bm{g}^{(i)}\left(\bm \uptheta\right) $ and $\hat{\bm{g}}_{i}\left(\bm \uptheta\right)$ in algorithms with \textit{distributed update} is the estimate of $\bm{g}_{i}\left(\bm \uptheta\right)$. If we choose the same SA (e.g., stochastic gradient) method to estimate the gradient in all the algorithms, the computational step is the same. In addition, the communication step is not necessarily required because for algorithms with \textit{cyclic update} because each agent only updates its particular components of the global decision variables and the next update does not depend on the information of other agents. As for algorithms with \textit{distributed update}, the gradient estimate only depends on local decision variables and whether communicating with other agents will not affect the updated results. In this case, algorithms with \textit{cyclic update} can also be implemented in parallel.


\section{Conclusion}
\label{sec:conclusion}
This paper compares the cyclic and the distributed approach in solving stochastic optimization problems with its main focus being on multi-agent systems. 
We show that the approximately
 are different from in both the framework of optimization problem and the 
 update rule. In particular, we use two examples to better exemplify different problems that the two approaches can solve. 

The main goal of the paper is to show the differences between the GCSA algorithm and the DSA algorithm. However, we present two other algorithms to help illustrate the distinctions of them. Clearly defining the \textit{cyclic framework}, \textit{distributed framework}, \textit{cyclic update}, and \textit{distributed update} helps researchers and practitioners avoid confusion when applying these algorithms to multi-agent problems. 

\bibliographystyle{plain}
\bibliography{sample}

\begin{thebibliography}{10}

\bibitem{beard2002coordinated}
Randal~W Beard, Timothy~W McLain, Michael~A Goodrich, and Erik~P Anderson.
\newblock Coordinated target assignment and intercept for unmanned air vehicles.
\newblock {\em IEEE Transactions on Robotics and Automation}, 18(6):911--922, 2002.

\bibitem{berahas2023balancing}
Albert~S Berahas, Raghu Bollapragada, and Shagun Gupta.
\newblock Balancing communication and computation in gradient tracking algorithms for decentralized optimization.
\newblock {\em arXiv preprint arXiv:2303.14289}, 2023.

\bibitem{berahas2018balancing}
Albert~S Berahas, Raghu Bollapragada, Nitish~Shirish Keskar, and Ermin Wei.
\newblock Balancing communication and computation in distributed optimization.
\newblock {\em IEEE Transactions on Automatic Control}, 64(8):3141--3155, 2018.

\bibitem{berahas2019derivative}
Albert~S Berahas, Richard~H Byrd, and Jorge Nocedal.
\newblock Derivative-free optimization of noisy functions via quasi-newton methods.
\newblock {\em SIAM Journal on Optimization}, 29(2):965--993, 2019.

\bibitem{bertsekasandtsitsiklis1989}
D.~P. Bertsekas and J.~N. Tsitsiklis.
\newblock {\em Parallel and Distributed Computation: Numerical methods}.
\newblock Prentice-Hall, 1989.

\bibitem{bianchi2013performance}
Pascal Bianchi, Gersende Fort, and Walid Hachem.
\newblock Performance of a distributed stochastic approximation algorithm.
\newblock {\em IEEE Transactions on Information Theory}, 59(11):7405--7418, 2013.

\bibitem{borkar1998asynchronous}
Vivek~S Borkar.
\newblock Asynchronous stochastic approximations.
\newblock {\em SIAM Journal on Control and Optimization}, 36(3):840--851, 1998.

\bibitem{botts2020three}
Carsten~H. Botts.
\newblock Three-dimensional swarming using cyclic stochastic optimization.
\newblock {\em IEEE Transactions on Aerospace and Electronic Systems}, 58(2):1431--1445, 2022.

\bibitem{botts2016multi}
Carsten~H Botts, James~C Spall, and Andrew~J Newman.
\newblock Multi-agent surveillance and tracking using cyclic stochastic gradient.
\newblock In {\em Proceedings of the 2016 American Control Conference (ACC)}, pages 270--275, 2016.

\bibitem{boyd2011distributed}
Stephen Boyd, Neal Parikh, and Eric Chu.
\newblock {\em Distributed optimization and statistical learning via the alternating direction method of multipliers}.
\newblock Now Publishers Inc, 2011.

\bibitem{bullo2018lectures}
Francesco Bullo, Jorge Cort{\'e}s, Florian D{\"o}rfler, and Sonia Mart{\'\i}nez.
\newblock {\em Lectures on network systems}, volume~1.
\newblock CreateSpace, 2018.

\bibitem{canutescu2003}
A.~A. Canutescu and R.~L. Dunbrack.
\newblock Cyclic coordinate descent: a robotics algorithm for protein loop closure.
\newblock {\em Protein Science}, 12(5):pp. 963--972, 2003.

\bibitem{cao2023first}
Liyuan Cao, Albert~S Berahas, and Katya Scheinberg.
\newblock First-and second-order high probability complexity bounds for trust-region methods with noisy oracles.
\newblock {\em Mathematical Programming}, pages 1--52, 2023.

\bibitem{cao2012overview}
Yongcan Cao, Wenwu Yu, Wei Ren, and Guanrong Chen.
\newblock An overview of recent progress in the study of distributed multi-agent coordination.
\newblock {\em IEEE Transactions on Industrial informatics}, 9(1):427--438, 2012.

\bibitem{chandler2001uav}
Phillip~R Chandler, Meir Pachter, and Steven Rasmussen.
\newblock {UAV} cooperative control.
\newblock In {\em Proceedings of the 2001 American Control Conference}, volume~1, pages 50--55, 2001.

\bibitem{chen2018stochastic}
Ruobing Chen, Matt Menickelly, and Katya Scheinberg.
\newblock Stochastic optimization using a trust-region method and random models.
\newblock {\em Mathematical Programming}, 169:447--487, 2018.

\bibitem{derakhshan2019review}
Farnaz Derakhshan and Shamim Yousefi.
\newblock A review on the applications of multiagent systems in wireless sensor networks.
\newblock {\em International Journal of Distributed Sensor Networks}, 15(5), 2019.

\bibitem{dong2018theory}
Xiwang Dong, Yongzhao Hua, Yan Zhou, Zhang Ren, and Yisheng Zhong.
\newblock Theory and experiment on formation-containment control of multiple multirotor unmanned aerial vehicle systems.
\newblock {\em IEEE Transactions on Automation Science and Engineering}, 16(1):229--240, 2018.

\bibitem{dorri2018multi}
Ali Dorri, Salil~S Kanhere, and Raja Jurdak.
\newblock Multi-agent systems: A survey.
\newblock {\em Ieee Access}, 6:28573--28593, 2018.

\bibitem{duchi2011adaptive}
John Duchi, Elad Hazan, and Yoram Singer.
\newblock Adaptive subgradient methods for online learning and stochastic optimization.
\newblock {\em Journal of machine learning research}, 12(7), 2011.

\bibitem{forero2010consensus}
Pedro~A Forero, Alfonso Cano, and Georgios~B Giannakis.
\newblock Consensus-based distributed linear support vector machines.
\newblock In {\em Proceedings of the 9th ACM/IEEE international conference on information processing in sensor networks}, pages 35--46, 2010.

\bibitem{gazi2011swarm}
Veysel Gazi and Kevin~M Passino.
\newblock {\em Swarm stability and optimization}.
\newblock Springer Science \& Business Media, 2011.

\bibitem{ghadimi2013stochastic}
Saeed Ghadimi and Guanghui Lan.
\newblock Stochastic first-and zeroth-order methods for nonconvex stochastic programming.
\newblock {\em SIAM journal on optimization}, 23(4):2341--2368, 2013.

\bibitem{granichin2018cyclic}
Oleg~N Granichin and Viktoriya~A Erofeeva.
\newblock Cyclic stochastic approximation with disturbance on input in the parameter tracking problem based on a multiagent algorithm.
\newblock {\em Automation and Remote Control}, 79(6):1013--1028, 2018.

\bibitem{haaland2010statistical}
Ben Haaland, Wanli Min, Peter~ZG Qian, and Yasuo Amemiya.
\newblock A statistical approach to thermal management of data centers under steady state and system perturbations.
\newblock {\em Journal of the American Statistical Association}, 105(491):1030--1041, 2010.

\bibitem{cuevas2017cyclic}
Karla Hernandez.
\newblock Cyclic stochastic optimization: Generalizations, convergence, and applications in multi-agent systems.
\newblock {\em arXiv preprint arXiv:1707.06700}, 2017.

\bibitem{hernandez2014cyclic}
Karla Hern{\'a}ndez and James~C Spall.
\newblock Cyclic stochastic optimization with noisy function measurements.
\newblock In {\em 2014 American Control Conference}, pages 5204--5209. IEEE, 2014.

\bibitem{johnson2013accelerating}
Rie Johnson and Tong Zhang.
\newblock Accelerating stochastic gradient descent using predictive variance reduction.
\newblock {\em Advances in neural information processing systems}, 26, 2013.

\bibitem{kiefer1952stochastic}
Jack Kiefer and Jacob Wolfowitz.
\newblock Stochastic estimation of the maximum of a regression function.
\newblock {\em The Annals of Mathematical Statistics}, pages 462--466, 1952.

\bibitem{kingma2014adam}
Diederik~P Kingma.
\newblock Adam: A method for stochastic optimization.
\newblock {\em arXiv preprint arXiv:1412.6980}, 2014.

\bibitem{kushner2012stochastic}
Harold~Joseph Kushner and Dean~S Clark.
\newblock {\em Stochastic approximation methods for constrained and unconstrained systems}, volume~26.
\newblock Springer Science \& Business Media, 2012.

\bibitem{larson2016stochastic}
Jeffrey Larson and Stephen~C Billups.
\newblock Stochastic derivative-free optimization using a trust region framework.
\newblock {\em Computational Optimization and applications}, 64:619--645, 2016.

\bibitem{larson2024novel}
Jeffrey Larson, Matt Menickelly, and Jiahao Shi.
\newblock A novel noise-aware classical optimizer for variational quantum algorithms.
\newblock {\em arXiv preprint arXiv:2401.10121}, 2024.

\bibitem{leenpark2008}
S.~Lee and F.~C. Park.
\newblock Cyclic optimization algorithms for simultaneous structure and motion recovery in computer vision.
\newblock {\em Engineering Optimization}, 40(5):pp. 403--419, 2008.

\bibitem{liandpetriouli2014}
B.~Li and A.~Petropuli.
\newblock Efficient target estimation in distributed {MIMO} radar via the {ADMM}.
\newblock In {\em Proceedings of the Conference on Information Sciences and Systems (CISS)}, pages 1--5, 2014.

\bibitem{liandosher2009}
Y.~Li and S.~Osher.
\newblock Coordinate descent optimization for $\ell_1$ minimization with application to compressed sensing; a greedy algorithm.
\newblock {\em Inverse Problems and Imaging}, 3(3):pp. 487--503, 2009.

\bibitem{liu2014asynchronous}
Ji~Liu, Steve Wright, Christopher R{\'e}, Victor Bittorf, and Srikrishna Sridhar.
\newblock An asynchronous parallel stochastic coordinate descent algorithm.
\newblock In {\em International Conference on Machine Learning}, pages 469--477. Proceedings of Machine Learning Research, 2014.

\bibitem{mcarthur2007multi}
Stephen~DJ McArthur, Euan~M Davidson, Victoria~M Catterson, Aris~L Dimeas, Nikos~D Hatziargyriou, Ferdinanda Ponci, and Toshihisa Funabashi.
\newblock Multi-agent systems for power engineering applications—part i: Concepts, approaches, and technical challenges.
\newblock {\em IEEE Transactions on Power systems}, 22(4):1743--1752, 2007.

\bibitem{mcarthur2007multi2}
Stephen~DJ McArthur, Euan~M Davidson, Victoria~M Catterson, Aris~L Dimeas, Nikos~D Hatziargyriou, Ferdinanda Ponci, and Toshihisa Funabashi.
\newblock Multi-agent systems for power engineering applications—part ii: Technologies, standards, and tools for building multi-agent systems.
\newblock {\em IEEE Transactions on Power Systems}, 22(4):1753--1759, 2007.

\bibitem{nedic2017achieving}
Angelia Nedic, Alex Olshevsky, and Wei Shi.
\newblock Achieving geometric convergence for distributed optimization over time-varying graphs.
\newblock {\em SIAM Journal on Optimization}, 27(4):2597--2633, 2017.

\bibitem{nedic2009distributed}
Angelia Nedic and Asuman Ozdaglar.
\newblock Distributed subgradient methods for multi-agent optimization.
\newblock {\em IEEE Transactions on Automatic Control}, 54(1):48--61, 2009.

\bibitem{nelder1965simplex}
John~A Nelder and Roger Mead.
\newblock A simplex method for function minimization.
\newblock {\em The computer journal}, 7(4):308--313, 1965.

\bibitem{nesterov2012efficiency}
Yu~Nesterov.
\newblock Efficiency of coordinate descent methods on huge-scale optimization problems.
\newblock {\em SIAM Journal on Optimization}, 22(2):341--362, 2012.

\bibitem{olfati2007consensus}
Reza Olfati-Saber, J~Alex Fax, and Richard~M Murray.
\newblock Consensus and cooperation in networked multi-agent systems.
\newblock {\em Proceedings of the IEEE}, 95(1):215--233, 2007.

\bibitem{ouyang2013stochastic}
Hua Ouyang, Niao He, Long Tran, and Alexander Gray.
\newblock Stochastic alternating direction method of multipliers.
\newblock In {\em International Conference on Machine Learning}, pages 80--88. Proceedings of Machine Learning Research, 2013.

\bibitem{peng2016arock}
Zhimin Peng, Yangyang Xu, Ming Yan, and Wotao Yin.
\newblock Arock: an algorithmic framework for asynchronous parallel coordinate updates.
\newblock {\em SIAM Journal on Scientific Computing}, 38(5):A2851--A2879, 2016.

\bibitem{peterson2014simulation}
Cameron~K Peterson, Andrew~J Newman, and James~C Spall.
\newblock Simulation-based examination of the limits of performance for decentralized multi-agent surveillance and tracking of undersea targets.
\newblock In {\em Signal Processing, Sensor/Information Fusion, and Target Recognition XXIII}, volume 9091, page 90910F. International Society for Optics and Photonics, 2014.

\bibitem{pipattanasomporn2009multi}
Manisa Pipattanasomporn, Hassan Feroze, and Saifur Rahman.
\newblock Multi-agent systems in a distributed smart grid: Design and implementation.
\newblock In {\em Proceedings of the IEEE/PES Power Systems Conference and Exposition}, pages 1--8, 2009.

\bibitem{powell2002uobyqa}
Michael~JD Powell.
\newblock {UOBYQA}: unconstrained optimization by quadratic approximation.
\newblock {\em Mathematical Programming}, 92(3):555--582, 2002.

\bibitem{pu2021distributed}
Shi Pu and Angelia Nedi{\'c}.
\newblock Distributed stochastic gradient tracking methods.
\newblock {\em Mathematical Programming}, 187(1):409--457, 2021.

\bibitem{ram2009incremental}
S~Sundhar Ram, A~Nedi{\'c}, and Venugopal~V Veeravalli.
\newblock Incremental stochastic subgradient algorithms for convex optimization.
\newblock {\em SIAM Journal on Optimization}, 20(2):691--717, 2009.

\bibitem{ram2010distributed}
S~Sundhar Ram, Angelia Nedi{\'c}, and Venugopal~V Veeravalli.
\newblock Distributed stochastic subgradient projection algorithms for convex optimization.
\newblock {\em Journal of Optimization Theory and Applications}, 147(3):516--545, 2010.

\bibitem{robbins1951stochastic}
Herbert Robbins and Sutton Monro.
\newblock A stochastic approximation method.
\newblock {\em The annals of mathematical statistics}, pages 400--407, 1951.

\bibitem{sahu2018distributed}
Anit~Kumar Sahu, Dusan Jakovetic, Dragana Bajovic, and Soummya Kar.
\newblock Distributed zeroth order optimization over random networks: A {K}iefer-{W}olfowitz stochastic approximation approach.
\newblock In {\em Proceeding of the 2018 IEEE Conference on Decision and Control (CDC)}, pages 4951--4958, 2018.

\bibitem{sergeenko2020advanced}
Anna Sergeenko, Oleg Granichin, and Anton~V Proskurnikov.
\newblock Advanced {SPSA}-based algorithm for multi-target tracking in distributed sensor networks.
\newblock In {\em Proceeding of the 2020 59th IEEE Conference on Decision and Control (CDC)}, pages 2424--2429, 2020.

\bibitem{shamma2007cooperative}
Jeff~S Shamma.
\newblock {\em Cooperative control of distributed multi-agent systems}.
\newblock Wiley Online Library, 2007.

\bibitem{shi2021sqp}
Jiahao Shi and James~C Spall.
\newblock Sqp-based projection spsa algorithm for stochastic optimization with inequality constraints.
\newblock In {\em 2021 American control conference (ACC)}, pages 1244--1249. IEEE, 2021.

\bibitem{shi2015extra}
Wei Shi, Qing Ling, Gang Wu, and Wotao Yin.
\newblock Extra: An exact first-order algorithm for decentralized consensus optimization.
\newblock {\em SIAM Journal on Optimization}, 25(2):944--966, 2015.

\bibitem{spallcyclicseesaw2012}
J.~C. Spall.
\newblock Cyclic seesaw process for optimization and identification.
\newblock {\em Journal of Optimization Theory and Applications}, 154(1):pp. 187--208, 2012.

\bibitem{spall1992multivariate}
James~C Spall.
\newblock Multivariate stochastic approximation using a simultaneous perturbation gradient approximation.
\newblock {\em IEEE transactions on automatic control}, 37(3):332--341, 1992.

\bibitem{spall2005introduction}
James~C Spall.
\newblock {\em Introduction to stochastic search and optimization: estimation, simulation, and control}, volume~65.
\newblock John Wiley \& Sons, 2005.

\bibitem{sun2023trust}
Shigeng Sun and Jorge Nocedal.
\newblock A trust region method for noisy unconstrained optimization.
\newblock {\em Mathematical Programming}, 202(1):445--472, 2023.

\bibitem{tseng2001}
P.~Tseng.
\newblock Convergence of a block coordinate descent method for non-differentiable minimization.
\newblock {\em Journal of Optimization Theory and Applications}, 109(3):pp. 475--494, 2001.

\bibitem{tsianos2012consensus}
Konstantinos~I Tsianos, Sean Lawlor, and Michael~G Rabbat.
\newblock Consensus-based distributed optimization: Practical issues and applications in large-scale machine learning.
\newblock In {\em Proceedings of the 50th Annual Allerton Conference on Communication, Control, and Computing}, pages 1543--1550. IEEE, 2012.

\bibitem{tsitsiklis1986distributed}
John Tsitsiklis, Dimitri Bertsekas, and Michael Athans.
\newblock Distributed asynchronous deterministic and stochastic gradient optimization algorithms.
\newblock {\em IEEE Transactions on Automatic Control}, 31(9):803--812, 1986.

\bibitem{tsitsiklis1984problems}
John~N Tsitsiklis.
\newblock {\em Problems in decentralized decision making and computation}.
\newblock PhD thesis, Massachusetts Institute of Technology, 1984.

\bibitem{tsitsiklis1994asynchronous}
John~N Tsitsiklis.
\newblock Asynchronous stochastic approximation and {Q}-learning.
\newblock {\em Machine learning}, 16(3):185--202, 1994.

\bibitem{wooldridge2009introduction}
Michael Wooldridge.
\newblock {\em An introduction to multiagent systems}.
\newblock John Wiley \& Sons, 2009.

\bibitem{wright2015}
S.~J. Wright.
\newblock Coordinate descent algorithms.
\newblock {\em Mathematical Programming}, 151(1):pp. 3--34, 2015.

\bibitem{yagan2006distributed}
Daniel Yagan and Chen-Khong Tham.
\newblock Distributed model-free stochastic optimization in wireless sensor networks.
\newblock In {\em International Conference on Distributed Computing in Sensor Systems}, pages 85--100. Springer, 2006.

\bibitem{yang2016parallel}
Yang Yang, Gesualdo Scutari, Daniel~P Palomar, and Marius Pesavento.
\newblock A parallel decomposition method for nonconvex stochastic multi-agent optimization problems.
\newblock {\em IEEE Transactions on Signal Processing}, 64(11):2949--2964, 2016.

\bibitem{zhu2020error}
Jingyi Zhu.
\newblock Error bounds and applications for stochastic approximation with non-decaying gain.
\newblock {\em arXiv preprint arXiv:2003.07357}, 2020.

\end{thebibliography}

\appendix

\end{document}